\documentclass{amsart}

\usepackage{amsmath,amssymb,amscd,amsfonts}

\newtheorem{theorem}{Theorem}
\newcommand{\bt}{\begin{theorem}}
\newcommand{\et}{\end{theorem}}
\newtheorem{lemma}{Lemma}
\newcommand{\bl}{\begin{lemma}}
\newcommand{\el}{\end{lemma}}
\newtheorem{corollary}{Corollary}
\newcommand{\bc}{\begin{corollary}}
\newcommand{\ec}{\end{corollary}}
\newcommand{\beq}{\begin{equation}}
\newcommand{\eeq}{\end{equation}}
\newcommand{\benum}{\begin{enumerate}}
\newcommand{\eenum}{\end{enumerate}}
\newcommand{\N}{\ensuremath{ \mathbf N }}
\newcommand{\Z}{\ensuremath{\mathbf Z}}
\newcommand{\Q}{\ensuremath{\mathbf Q}}
\newcommand{\R}{\ensuremath{\mathbf R}}

\newcommand{\mcp}{\ensuremath{ \mathcal P}}

\DeclareMathOperator{\id}{id}

\begin{document}

\title{Cantor polynomials for semigroup sectors }
\author{Melvyn B. Nathanson}
\address{Department of Mathematics\\
Lehman College (CUNY)\\
Bronx, NY 10468}
\email{melvyn.nathanson@lehman.cuny.edu}


\date{\today}

\begin{abstract}
A packing function on a set $\Omega$ in $\R^n$  is a one-to-one 
correspondence between  the set of lattice points in $\Omega$ 
and the set $\N_0$ of nonnegative integers.  
It is proved that if $r$ and $s$ are relatively prime positive integers 
such that $r$ divides $s-1$, then there exist two distinct quadratic packing 
polynomials on the sector $\{ (x,y) \in \R^2 : 0 \leq y \leq rx/s\}$. 
For the rational numbers $1/s$, these are the unique quadratic packing polynomials. 
Moreover, quadratic  quasi-polynomial packing functions 
are constructed for all rational sectors.  
\end{abstract}

\maketitle

\section{The Cantor polynomial packing problem}

A \emph{packing function} 
\index{packing function} on a set $\Omega$ in $\R^n$  is a function
 that is a one-to-one 
correspondence between  the set of lattice points in $\Omega$ 
and the set $\N_0$ of nonnegative integers.  
A \emph{packing polynomial} \index{packing polynomial} on $\Omega$ 
is a packing function that is a polynomial.  
Packing polynomials are  used in computer science to store and retrieve 
multi-dimensional arrays in linear memory. 
In this paper we consider only $n = 2$.  

Define $e_{2,s} = (s,1)$  for all $s \in \N_0$.
Let $e_1 = (1,0)$ and $e_2 = e_{2,0} = (0,1)$ 
The \emph{standard basis} for $\R^2$ is $\{e_1,e_2\}$.  
We define $1/0 = \infty$ and $1/\infty = 0$. 
 
For every positive real number $\alpha$ we consider the 
\emph{real sector}  \index{real sector}
\[
S(\alpha) = \{ (x,y) \in \R^2 : 0 \leq y \leq \alpha x\}
\]
and the \emph{integer sector}  \index{integer sector}
\[
I(\alpha) = S(\alpha) \cap \N_0^2 
= \{ (x,y) \in \N_0^2 : 0 \leq y \leq \alpha x\}.
\]
If $(a,b) \in S(\alpha)\setminus \{ (0,0)\}$, then $b/a \leq \alpha$.
A sector is called \emph{rational}  \index{rational sector} if $\alpha$ is a rational number 
and \emph{irrational}  \index{irrational sector} if $\alpha$ is an irrational number.  
The real sector $S(\alpha)$ is the cone with vertex at $(0,0)$ 
generated by the points $(1,0)$ and $(1,\alpha)$.  Equivalently, $S(\alpha)$ is 
the convex hull of the rays $\{ (x,0) : x \geq 0\}$ and $\{(x,\alpha x) : x \geq 0\}$.

We also define the real and integer sectors
\[
S(\infty) = \{ (x,y) \in \R^2 : x \geq 0 \text{ and } y \geq 0\}
\]
and
\[
I(\infty) = S(\infty) \cap \N_0^2 = \N_0^2.
\]
The real sector $S(\infty)$ is the cone with vertex at $(0,0)$ 
generated by the points $(1,0)$ and $(0,1)$.

Cantor proved that the sector $I(\infty)$ has the same 
cardinality as the set $\N_0$  by constructing two explicit 
packing polynomials from $I(\infty)$ to $\N_0$:
\[
F_{\infty}(x,y) = \frac{(x+y)^2}{2}+ \frac{x+3y}{2}
\]
and 
\[
G_{\infty}(x,y) = \frac{(x+y)^2}{2}+ \frac{3x+y}{2}
\]
These are called the \index{Cantor polynomials} \emph{Cantor  polynomials}.
Using methods from  analytic number theory, 
Fueter and P\' olya~\cite{fuet-poly23} proved in 1923 that these 
are the unique quadratic 
polynomials that are bijections from  $I(\infty)$ to $\N_0$.  
Recently, Vsemirnov~\cite{vsem01} gave an elementary proof 
using only quadratic reciprocity.  
The survey paper~\cite{nath13xy} contains an exposition of Vsemirnov's proof and other results on packing polynomials.  

It has been conjectured since the work of Fueter and P\' olya that the two 
Cantor polynomials are the only  packing polynomials on $I(\infty)$, or, 
equivalently, that every packing polynomial on $I(\infty)$ is quadratic.  
Lew and Rosenberg~\cite{lew-rose78b} proved that there is no cubic or quartic 
packing polynomial on $I(\infty)$.
 Despite considerable work  
(see references), the Cantor polynomial problem is still unsolved.  
In this paper we construct two quadratic packing polynomials on $I(\alpha)$ 
for certain rational numbers $\alpha$, and apply the Fueter-P\' olya theorem 
to prove that these are the only quadratic packing polynomials on these sectors.  
We also construct quadratic  quasi-polynomial packing functions 
for all rational sectors.  
The final section lists some related open problems.

\section{Sectors and semigroups}
A subset $W$ of an additive abelian semigroup is a \emph{generating set} for the 
semigroup if every element of the semigroup can be represented  
as a finite  nonnegative integral linear combination of elements of $W$.  
Such additive representations are not necessarily unique.  
An additive abelian  semigroup is  \emph{free of rank $k$}  
if it contains a generating set $W = \{w_1,\ldots, w_k\}$ such that 
every element in the semigroup has a unique representation 
in the form $\sum_{i=1}^k x_i w_i$ with $x_i \in \N_0$ 
for $i=1,\ldots, k$.  
The set $W = \{w_1,\ldots, w_k\}$ is called a \emph{free basis} for the semigroup.
If $w$ is an element  of a free basis, then $w \neq 0$.

For every real number $\alpha > 0$ and for $\alpha = \infty$, 
the real sector $S(\alpha)$ and the integer sector $I(\alpha)$ 
are additive abelian semigroups.  The following lemma 
identifies the numbers $\alpha$ such that $I(\alpha)$ is  a free semigroup.

\bt      \label{CantorSector:theorem:FreeBasis}
The additive abelian semigroup $I(\alpha)$ is free of rank $k$ 
if and only if $k=2$ and $\alpha \in \{ 1/s: s\in \N_0\}$. 
Moreover,  $\{(1,0), (s,1) \}$ is the unique free basis of $I(1/s)$ for all $s \in \N_0$.
\et

\begin{proof}
Let $\alpha$ be a positive real irrational number, 
and let $W$ be a finite subset of $I(\alpha)\setminus \{ (0,0) \}$.
Let $W = \{ w_1,\ldots, w_k\}$, where 
$w_i = (a_i,b_i) \in I(\alpha)\setminus \{ (0,0) \}$ and $a_i \neq 0$  for $i=1,\ldots, k$.  
We define 
\[
\lambda = \min\left\{ \frac{b_i}{a_i} : i = 1, \ldots, k  \right\}
\]
and
\[
\mu = \max\left\{ \frac{b_i}{a_i} : i = 1, \ldots, k  \right\}.
\]
The cone generated by the nonnegative rays $y = \lambda x$ and $y = \mu x$ is
\[
C =  \{ (x,y) \in S(\infty) : \lambda x \leq y \leq \mu x\}.
\]
This cone contains $W$, and so contains the additive subgroup generated by $W$.

Because $b_i/a_i \leq \alpha$ and $\alpha$ is irrational, 
we have $b_i/a_i < \alpha$ for all $i=1,\ldots, k$, 
and so $\mu < \alpha$.
There exist positive integers $c$ and $d$ such that 
\[
\mu < \frac{c}{d} < \alpha
\]
and the lattice point $(c,d)$ belongs to $I(\alpha)$ but not to $C$. 
Thus, the semigroup $I(\alpha)$ is not finitely generated.  
In particular, $I(\alpha)$ is not a free abelian semigroup of rank $k$ for any $k \in \N$.

Let $\Gamma \neq \{(0,0)\}$ be an additive abelian semigroup 
contained in $\N_0^2$ or, more generally, 
in $\Q_0^2$, where $\Q_0$ is the set of nonnegative rational numbers.  
If $\Gamma$ is free of rank $k \geq 3$, then there exist $w_1,w_2,w_3 \in \Gamma$ 
that are elements of a free basis for $\Gamma$.
Because 
\[
w_i \in \Q_0^2 \setminus \{ (0,0) \} \subseteq \Q^2 \setminus \{ (0,0) \}
\]
 for $i=1,2,3$, it follows that $w_1,w_2,w_3$ 
are \Q-linearly dependent, and there exist rational numbers $t_1,t_2,t_3$, 
not all 0, such that $t_1w_1+t_2w_2+t_3w_3 = (0,0).$  
The vectors $w_i$ have nonnegative coordinates, and so 
$t_i > 0$ for some $i$ and $t_j < 0$ for some $j$.  
Thus, we can assume that $t_1>0, t_2 \geq 0$, and $t_3 < 0$.  We have  
\[
t_1w_1+t_2w_2 =-t_3w_3.
\] 
Multiplying by a common multiple of the denominators of the fractions $t_1,t_2,t_3$, 
we obtain nonnegative integers $x_1,x_2,x_3$ with $x_1 \geq 1$ and $x_3 \geq 1$ 
such that 
\[
x_1w_1+x_2w_2 = x_3w_3.
\]
This is impossible if $w_1,w_2,w_3$ belong to a free basis for $\Gamma$.
Thus, $k \leq 2$.  If $k=0$, then $\Gamma = \{ (0,0)\}$, which is absurd. 
Thus, if $\Gamma$ is a free abelian semigroup of rank $k$, then $k=1$ or 2.

If $\alpha > 0$ or $\alpha = \infty$, then 
$(x,0) \in I(\alpha)$ for all $x \in \N_0$, and $(x,1) \in I(\alpha)$ for all $x \geq 1/\alpha$.  
Thus, $I(\alpha)$ is an abelian semigroup in $\Q_0^2$ 
that does not lie on a line, and so, if $I(\alpha)$ is free of rank $k$, then $k=2$.

Let $r$ and $s$ be relatively prime positive integers with $r \geq 2$.  
The semigroup $I(r/s)$ is contained in the convex hull of the nonnegative rays $y = 0$ 
and $y = rx/s$.  
The lattice points $(1,0)$  and $(s,r)$ are in $I(r/s)$ 
and lie on the rays $y=0$ and $y  = rx/s$, respectively.  
Let $W = \{w_1,w_2\} \subseteq I(r/s)\setminus \{(0,0)\}$.  
If $W$ generates $I(r/s)$, then there exist nonnegative integers $x_1$ 
and $x_2$ such that  such that $(1,0) = x_1w_1 + x_2w_2$ for $x_1,x_2\in \N_0$.
This is possible only if $w_i = (1,0)$ for $i = 1$ or 2.  Let $w_1 = (1,0)$.  
Because $W$ generates $I(r/s)$, then there exist nonnegative integers 
$y_1$ and $y_2$ such that  such that $(s,r) = y_1(1,0) + y_2w_2$.  
It follows that $y_2w_2 = (s,r) - y_1(1,0) =  (s - y_1,r) \in I(r/s)$.
However, if $(a,r) \in I(r/s)$, then $a \geq s$ and so $y_1=0$ and 
$(s,r) = y_2w_2$.  Because $r$ and $s$ are relatively prime, 
we conclude that $y_2 = 1$ and $w_2 = (s,r)$.  Thus, 
every element generated by $W$ is of the form 
\[
x_1(1,0) + x_2(s,r) = (x_1+x_2s, x_2r)
\]
with $x_1, x_2 \in \N_0$.  
Because $(x,1) \in I(r/s)$ for all $x \geq s/r$, it follows that $x_2r=1$ for some
$x_2 \in \N_0$, and this is possible only if $x_2 = r = 1$.
Thus, if the integer sector $I(r/s)$ is generated by a set $W$ of cardinality 2, 
then $r=1$ and $W = \{ (1,0), (s,1) \}$.  

If $(x,y) \in I(1/s)$, then $y \leq x/s$ or, equivalently, $x - sy \in \N_0$, and 
\[
(x,y) = ((x-sy)+sy,y) = (x-sy,0)+ (sy,y) = (x-sy)(1,0) + y(s,1).
\]
This representation of $(x,y)$ is unique because the vectors $(1,0)$ and 
$(s,1)$ are linearly independent in $\R^2$.
It follows that $I(1/s)$ is free of rank 2, and that $\{ (1,0), (s,1) \}$ 
is the unique basis for $I(1/s)$.  

Finally, the semigroup $I(\infty) = \N_0^2$ is free, and $\{ (1,0), (0,1)\}$ 
is the unique free basis for $I(\infty)$.  
This completes the proof.  
\end{proof}

\bt       \label{CantorSector:theorem:FreeBasisTrans}
There exists a linear transformation from $\R^2$ to $\R^2$ 
whose restriction to $I(\infty)$ is a one-to-one function onto $I(\alpha)$ 
if and only if $\alpha \in \{ 1/s: s\in \N_0\}$.  
Moreover, for each $s \in \N_0$  there are exactly two such linear transformations, 
whose matrices with respect to the standard basis are 
\beq     \label{CantorSector:LambdaMs}
\Lambda_s = \left(
\begin{matrix}
1 & s \\
0 & 1
\end{matrix}
\right)
\qquad \text{and} \qquad
M_s = \left(
\begin{matrix}
s & 1 \\
1 & 0
\end{matrix}
\right).
\eeq
\et

\begin{proof} 
Let $T:\R^2 \rightarrow \R^2$ be a linear transformation.  
Let $e_1 = (1,0)$ and $e_2 = (0,1)$, and let $w_1 = T(e_1)$ and $w_2 = T(e_2)$.  
Then $I(\infty)) = \{xe_1+ye_2 : x,y \in \N_0\}$ and 
\[
T(I(\infty)) = \{xT(e_1)+yT(e_2) : x,y \in \N_0 \} 
= \{xw_1+yw_2 : x,y \in \N_0 \}. 
\]
Thus, $T(I(\infty))$ is the additive abelian semigroup generated by the set $\{ w_1, w_2\}$.
If the restriction of $T$ to $I(\infty)$ is a one-to-one function, 
then $T(I(\infty))$ is a free semigroup of rank 2, 
and $\{ w_1, w_2\}$ is a free basis for $T(I(\infty)) $.  
If  the restriction of $T$ to $I(\infty)$ is a one-to-one function onto $I(\alpha)$, 
then $I(\alpha)$ is a free semigroup of rank 2, 
and  $\alpha \in \{ 1/s: s\in \N_0 \}$ 
by Theorem~\ref{CantorSector:theorem:FreeBasis}.

The unique free basis for $I(\infty)$ is $\{ e_1,e_2\}$.  
If  the restriction of $T$ to $I(\infty)$ is a one-to-one function onto $I(\infty)$, 
then $T = \Lambda_0 = \id$ or $T = M_0$.  
For every positive integer $s$, the unique free basis for $I(1/s)$ is $\{ e_1,e_{2,s}\}$.  
If the restriction of $T$ to $I(\infty)$ is 
a one-to-one function onto $I(1/s)$,  then $T = \Lambda_s$ or $T = M_s$.  
This completes the proof.
\end{proof}

For all nonnegative integers $s$ and $t$, we have 
$\Lambda_{s+t} = \Lambda_s\Lambda_t$, and so 
$\Lambda_s = \Lambda_1^s$ and  $\Lambda_s^{-1} = \Lambda_{-s}$.
Moreover,  $\Lambda_s$ is a bijection from $I(1/t)$ to $I(1/(s+t) )$, 
and $M_s = \Lambda_s M_0.$

\bt                            \label{CantorSector:theorem:FreeBasisInvolut}
For every nonnegative integer $s$, the linear transformation 
\beq     \label{CantorSector:PhiInvolution}
\Phi_s = \left(
\begin{matrix}
s & 1-s^2 \\
1 & -s
\end{matrix}
\right)
\eeq
is the unique non-identity  linear transformation whose restriction to $I(1/s)$ 
is a bijection onto $I(1/s)$.  Moreover, $\Phi_s$ is an involution.  

For every positive integer $r$, the linear transformation 
\beq     \label{CantorSector:PsiInvolution}
\Psi_r = \left(
\begin{matrix}
1 & 0 \\
r & -1
\end{matrix}
\right)
\eeq
is an involution whose restriction to $I(r)$ 
is a bijection onto $I(r)$.  
\et

Note that $\Phi_1 = \Psi_1 = \left(
\begin{matrix}
1 & 0 \\
1 & -1
\end{matrix}
\right).$

\begin{proof}
By Theorem~\ref{CantorSector:theorem:FreeBasisTrans}, 
the unique free basis for $I(1/s)$ is $\{ e_1,e_{2,s} \}$.  
If $\Phi_s$ is a non-identity linear transformation that induces 
a bijection on $I(1/s)$, then $\Phi_s(e_1) = e_{2,s}$ 
and $\Phi_s(e_{2,s}) = e_1$.  Because $e_2  = e_{2,s} - se_1$, 
we have 
\[
\Phi_s(e_2) = \Phi_s\left(  e_{2,s} \right)- s\Phi_s\left( e_1 \right) 
= (1,0) - s(s,1) = (1-s^2, -s).
\]
This determines the matrix for $\Phi_s$ with respect to the standard basis.
Moreover, $\Phi_s^2$ is the identity transformation.

Similarly, $\Psi_r^2 = I$ and $\Psi_r(x,y) = (x , rx-y)$.  
If $0 \leq y \leq rx$, then $0 \leq rx-y \leq rx$, and so the 
restriction of $\Psi_r$ to $I(r)$ is a bijection.  
This completes the proof.  
\end{proof}

\section{Packing polynomials on rational sectors}

Let $\alpha > 0$ or $\alpha = \infty$.  
For every positive integer $d$, 
let $\mcp_d(\alpha)$ denote the set of packing polynomials 
of degree $d$ on the integer sector $I(\alpha)$.  
The Fueter-P\' olya theorem states that $\mcp(\infty) = \{F_{\infty}, G_{\infty}\}$. 
We also have $\mcp_3(\infty) = \mcp_4(\infty) =\emptyset$ by the Lew-Rosenberg theorem.

\bt
For $\alpha > 0$ and $\alpha = \infty$, 
there exists no linear packing polynomial on $I(\alpha)$.  
Equivalently,  $\mcp_1(\alpha) = \emptyset$.  
\et

\begin{proof}
This is a simple counting argument.  Let $\alpha > 0$.  
For every positive integer $n$, let 
\[
I_n(\alpha) = \{ (x,y) \in I(\alpha)  : x \leq n \}.
\]
Then
\[
|I_n(\alpha)| = \sum_{j=0}^n \left( \left[ \alpha j \right] + 1\right)
> \sum_{j=0}^n  \alpha j = \frac{\alpha n(n+1)}{2} > \frac{\alpha n^2}{2}.
\]
Let $f(x,y) = ax+by+c$ be a linear polynomial 
such that $f(x,y) \in \N_0$ for all $(x,y) \in I(\alpha)$.
If $(x,y) \in I_n(\alpha)$, then $x \leq n$ and $y \leq \alpha n$,  and so 
\[
0 \leq f(x,y) = ax+by+c 
\leq  \left(  |a| +  \alpha  \left| b\right|  + |c| \right) n.  
\]
Thus, the linear function $f$ maps a set with more than $(\alpha/2) n^2$ elements
into a set with at most  $C n$ elements, 
where $C =  |a| +  \alpha  \left| b\right|  + |c| + 1 $,  
and this function cannot be one-to-one for $n > 2C/\alpha$.  
It follows that there is no linear packing polynomial on $I(\alpha)$. 

There is a similar counting argument for $I(\infty)$ by considering the finite set
$I_n(\infty) = \{ (x,y) \in I(\infty)  : x \leq n \text{ and } y \leq n\}$.
\end{proof}

\bt   \label{CantorSector:theorem:PolynomialBijection-1}
For $s\in \N_0$ , let $\Lambda_s$ and $M_s$ be the matrices defined by~\eqref{CantorSector:LambdaMs}.  
The functions from $\mcp_d(1/s)$ to $\mcp_d(0)$ defined by 
$F \mapsto F\circ \Lambda_s$ 
and $F \mapsto F\circ M_s$ are bijections with inverses 
$F \mapsto F\circ \Lambda_s^{-1}$ 
and $F \mapsto F\circ M_s^{-1}$, respectively.  
\et

\begin{proof}
This follows immediately from Theorem~\ref{CantorSector:theorem:FreeBasisTrans}, 
and the observation that if $F= F(x,y)$ is a polynomial of degree $d$, 
then $F\circ \Lambda_s$ and $F \circ M_s$ are polynomials of degree $d$.
\end{proof}

\bt   \label{CantorSector:theorem:PolynomialBijection-2}
For $s\in \N_0$, let $\Phi_s$ be the matrix defined 
by~\eqref{CantorSector:PhiInvolution}.  
The function from $\mcp_d(1/s)$ to $\mcp_d(1/s)$ defined by 
$F \mapsto F\circ \Phi_s$ is an involution.
\et

\begin{proof}
This follows immediately from Theorem~\ref{CantorSector:theorem:FreeBasisInvolut}.
\end{proof}

\bt    \label{CantorSector:theorem:I(r)}
Let $r$ be a positive integer.  The polynomials 
\[
F_r(x,y) = \frac{rx(x-1)}{2} + x + y
\]
and
\[
G_r(x,y) = \frac{rx(x+1)}{2} + x  - y.
\]
are packing polynomials on the integer sector $I(r) = \{(x,y)\in \N_0^2:0 \leq y \leq rx\}$.
\et

\begin{proof}
Consider the finite set $J_a = \{(a,b)\in \N_0^2:0 \leq b \leq ra\}$.
We have $|J_a| = ra+1$, and 
\[
I(r) = \bigcup_{a\in \N_0} J_a.
\]
Moreover, if $a, a' \in \N_0$ and $a\neq a'$, then $J_a \cap J_{a'} = \emptyset$.  

We construct the function $F_r(x,y)$ by enumerating the lattice points in $I(r)$ 
as follows:  For $a \geq 1$, count the elements of $J_a$ after the elements of $J_{a-1}$, 
and count the elements of $J_a$ from bottom $(a,0)$ to top $(a,ra)$.
Thus, $F_r(0,0) = 0$ and, for $x \geq 1$, 
\[
F_r(x,0) = 1 + \sum_{a=1}^{x-1} (ra+1) = \frac{rx(x-1)}{2} + x. 
\]
For $1 \leq y \leq rx$, 
\[
F_r(x,y) =  F_r(x,0)  + y = \frac{rx(x-1)}{2} + x + y.
\]

We construct the function $G_r(x,y)$ by enumerating the lattice points in $I(r)$ 
as follows:  For $a \geq 1$, count the elements of $J_a$ after the elements of $J_{a-1}$, 
and count the elements of $J_a$ from  top $(a,ra)$ to bottom $(a,0)$.
Thus, $G_r(0,0) = 0$ and, for $a \geq 1$, 
\[
G_r(x,rx) = 1 + \sum_{a=1}^{x-1} (ra+1) =  \frac{rx(x-1)}{2} + x.
\]
For $1 \leq y \leq x$, 
\[
G_r(x,y) = G_r(x,rx) + rx - y 
= \frac{rx(x+1)}{2} + x  - y.
\]
Thus, $F_r(x,y)$ and $G_r(x,y)$ are quadratic packing polynomials on $I(r)$.  
This completes the proof.  
\end{proof}

There is another way to construct the polynomial $G_r$.  
By Theorem~\ref{CantorSector:theorem:FreeBasisInvolut}, 
the linear transformation $\Psi_r:\R^2 \rightarrow \R^2$ that sends $(x,y)$ to 
$(x,rx-y)$ is a bijection on $I(r)$.  
Composing $F_r$ with $\Psi_r$, we obtain
\[
F_r\circ \Psi_r (x,y) = F_r(x,rx-y) = \frac{rx(x-1)}{2} + x + (rx-y) = G_r(x,y).
\]

\bt
$F_1(x,y)$ and $G_1(x,y)$ are the 
unique quadratic packing polynomials on $I(1)$.
\et

\begin{proof}
Applying Theorem~\ref{CantorSector:theorem:I(r)} with $r=1$, 
we obtain the quadratic packing polynomials   
\[
F_1(x,y) = \frac{x(x+1)}{2} +  y
\]
and
\[
G_1(x,y) = \frac{x(x+1)}{2} + x  - y
\]
on the integer sector $I(1)$.  
We also have
\[
F_{\infty}\circ \Lambda_1^{-1} (x,y) = F_{\infty}(x-y,y) = \frac{x^2}{2}+\frac{x+2y}{2} 
= F_1(x,y)
\]
and
\[
G_{\infty}\circ \Lambda_1^{-1} (x,y) = G_{\infty}(x-y,y) = \frac{x^2}{2}+\frac{3(x-y)+y}{2} 
= G_1(x,y).
\]
By Theorem~\ref{CantorSector:theorem:PolynomialBijection-1}, 
composition with  $\Lambda_1^{-1}$ is a bijection from 
$\mcp_2(0)$ to $\mcp_2(1)$.  
By the Fueter-P\' olya theorem, $\mcp_2(0) = \{F_{\infty}, G_{\infty} \}$ 
and so 
\[
\mcp_2(1) = \{ F_{\infty}\circ \Lambda_1^{-1} , G_{\infty}\circ \Lambda_1^{-1}  \} 
= \{F_1, G_1\}. 
\]
This completes the proof.  
\end{proof}

\bt    \label{CantorSector:theorem:IrsJ}
Let $r$ and $s$ be relatively prime positive integers 
such that $1 \leq r < s$ and $r$ divides $s-1$.  
Let $d = (s-1)/r$.  
For $a \in \N_0$, let
\[
J_a =  \left\{ ( a+dj, j ): j=0,1,\ldots, ra \right\}.
\]
Then
\[
I\left( r/s \right) 
= \bigcup_{a\in \N_0} J_a.
\]
Moreover, if $a \neq a'$, then $J_a \cap J_{a'} = \emptyset$.  
\et

\begin{proof}
If
\[
f_a(x) = \frac{x-a}{d}
\]
then $f_a(a+dj) = j $ and 
\[
J_a =  \left\{ \left( a+dj, f_a(a+dj) \right): j=0,1,\ldots, ra \right\}.
\]
If $a \neq a'$, the graphs of $y = f_a(x)$ and $y = f_{a'}(x)$ are distinct parallel lines, 
and so $J_a \cap J_{a'} = \emptyset$.

If  $(x,y) \in J_a$, then $x=a+dj$ and $y = j$ for some $j \in \{0,1,\ldots, ra\}$.
In particular, if $j = 0$, then $(x,y) = (a,0) \in S(r/s)$ and if $j = ra$, 
then $(x,y) = (a + dra,ra) = (sa,ra) \in S(r/s).$
It follows by convexity that the line segment between $(a,0)$ and $(sa,ra)$ 
lies in $S(r/s)$,
and this line segment is the graph of $y = f_a(x)$ for $a \leq x \leq sa$.
We conclude that $J_a \subseteq I(r/s)$ for all $a\in \N_0$, and so 
$\bigcup_{a\in \N_0} J_a \subseteq I(r/s)$.

Conversely, let $(x,y) \in I(r/s)$.  If $x=0$, then $y=0$ and $(x,y) = (0,0) \in J_0$.
If $x \geq 1$, then 
\[
y \leq \frac{r x}{s} < \frac{rx}{s-1} = \frac{x}{d} 
\]
and so 
\[
a = x - dy
\]
is a positive integer.  We have
\[
(x,y) = (a+dy,y).
\]
Moreover, $sy \leq rx$ implies that 
\[
y \leq rx - (s-1)y = r(x - dy ) = ra
\]
and so $(x,y) \in J_a.$  
It follows that $I(r/s)  \subseteq  \bigcup_{a\in \N_0} J_a$.  
This completes the proof.  
\end{proof}

\bt     \label{CantorSector:theorem:I(r/s)}
Let $r$ and $s$ be relatively prime positive integers 
such that $1 \leq r < s$ and $r$ divides $s-1$.  
Let $d = (s-1)/r$. 
The polynomials
\[
F_{r/s}(x,y) = \frac{r(x-dy)^2}{2} + \frac{ (2-r)x + (dr-2d+2)y }{2}
\]
and
\[
G_{r/s}(x,y) = \frac{r(x-dy)^2}{2} +  \frac{(r+2)x  - (2d+s+1)y}{2}.
\]
are quadratic packing polynomials for the sector $I(r/s)$.  
\et

\begin{proof}
Using Theorem~\ref{CantorSector:theorem:IrsJ}, we have
\[
|J_a| = ra+1
\]
and 
\[
1+\bigcup_{a=1}^{x-1} |J_a| = 1+ \sum_{a=1}^{x-1} (ra+1) = 
\frac{rx(x-1)}{2}+x.
\]
We construct the function $F_{r/s}(x,y)$ by enumerating the lattice points in $I(r/s)$ 
as follows:  For $a \geq 1$, count the elements of $J_a$ after the elements of $J_{a-1}$, 
and count the elements of $J_a$ from the bottom lattice point $(a,0)$ 
to the top lattice point $(ra,sa)$.
Thus, $F_{r/s}(0,0) = 0$ and, for $x \geq 1$, 
\[
F_{r/s}(x,0) = 1 + \sum_{a=1}^{x-1} (ra+1) = \frac{rx(x-1)}{2} + x.
\]
If $(x,y) \in I(r/s) \setminus \{ (0,0) \}$, then $0 \leq y \leq rx/s$ 
and $(x,y)=((x-dy)+dy,y)= (a+dy,y)$, 
where, as in the proof of Theorem~\ref{CantorSector:theorem:IrsJ}, the integer 
 $a=x-dy$ is positive, and so $(x,y) \in J_a$.  We have 
\begin{align*}
F_{r/s}(x,y) & = F_{r/s}(x-dy,0)+y \\ 
& = \frac{r(x-dy)(x-dy-1)}{2} + x- (d-1)y \\
& = \frac{r(x-dy)^2}{2} - \frac{r(x-dy)}{2} + x- (d-1)y.
\end{align*}

We construct the function $G_{r/s}(x,y)$ by enumerating the lattice points in $I(r/s)$ 
as follows:  For $a \geq 1$, count the elements of $J_a$ after the elements of $J_{a-1}$, 
and count the elements of $J_a$ from the top lattice point $(ra,sa)$ to 
the bottom lattice point $(a,0)$.
Thus, $G_{r/s}(0,0) = 0$ and, for $x \geq 1$, 
\[
G_{r/s}(sx,rx) = 1 + \sum_{a=1}^{x-1} (ra+1) = \frac{rx(x-1)}{2} + x.
\]
If $(x,y) \in I(r/s)\setminus \{ (0,0) \}$, then $0 \leq y \leq rx/s$ and $a = x-dy$ is a positive integer.  We have 
\[
(x,y)=((x-dy)+dy,y)= (a+dy,y) = (a+dy,ra - z)
\]
with $z = ra -y$.  
Then 
\begin{align*}
G_{r/s}(x,y) & = G_{r/s} (a+dy,ra - z) \\
& = G_{r/s}(sa,ra)+ z  \\ 
& =  \frac{ra(a-1)}{2} + a + z \\
& = \frac{r(x-dy)(x-dy-1)}{2} +(x-dy) + r(x-dy) - y \\
& = \frac{r(x-dy)(x-dy-1)}{2} +(r+1)x - (d + s)y \\
& = \frac{r(x-dy)^2}{2} + \frac{(r+2)x}{2}  - \frac{(2d+s+1)y}{2}.
\end{align*}
This completes the proof.  
\end{proof}

\bt     \label{CantorSector:theorem:I(1/s)}
For every integer $s \geq 2$, the polynomials
\[
F_{1/s}(x,y) = \frac{(x-(s-1)y)^2}{2} + \frac{x+(3-s)y}{2}
\]
and
\[
G_{1/s}(x,y) = \frac{(x-(s-1)y)^2}{2} + \frac{3x+(1-3s)y}{2}
\]
are the unique quadratic packing polynomials on the sector $I(1/s)$.  
\et

\begin{proof}
Applying Theorem~\ref{CantorSector:theorem:I(r/s)} with $r=1$ and $d = s-1$, 
we obtain the polynomials $F_{1/s}(x,y)$ and $G_{1/s}(x,y)$.  
We also have 
\begin{align*}
F_{\infty} \circ \Lambda_s^{-1}(x,y) 
& = F_{\infty}(x- sy,y) \\
& = \frac{(x-(s-1)y)^2}{2} +  \frac{ (x-sy)+3y }{2} \\
& = F_{1/s}(x,y).
\end{align*}
and
\begin{align*}
G_{\infty} \circ \Lambda_s^{-1}(x,y) 
& = G_{\infty}(x- sy,y) \\
& = \frac{(x-(s-1)y)^2}{2} +  \frac{ 3(x-sy)+y }{2}   \\
& = G_{1/s}(x,y).
\end{align*}
By Theorem~\ref{CantorSector:theorem:PolynomialBijection-1}, 
the function $F \mapsto F\circ \Lambda_s^{-1}$ is a bijection from 
$\mcp_2(0)$ onto $\mcp_2(1/s)$, and so 
 $\mcp_2(1/s) = \{F_{1/3}, G_{1/3}\}$.
Thus, $F_{1/3}$ and $G_{1/3}$ are the unique quadratic packing polynomials 
on $I(1/s)$.  
This completes the proof.  
\end{proof}

\section{Quasi-polynomials}
Let $H$ be a function whose domain is a set of integers.  
The function $H$ is a \emph{quasi-polynomial} with period $m$ if 
\[
H(x) = \sum_{i=0}^n c_i(x) x^i
\]
where the value of $c_i(x)$ depends only on the congruence class 
of $x$ modulo $m$.  
The coefficient functions $c_i(x)$ are not necessarily integer-valued.  
A quasi-polynomial restricted to a congruence class modulo $m$ is simply a polynomial.  
Every polynomial is a quasi-polynomial with period $m$ for every positive integer $m$.  

Similarly, if $H$ is a function whose domain is a set of lattice points in $\Z^2$,
then $H$ is a quasi-polynomial with period $m$ if
\[
H(x,y) = \sum_{i=0}^{d_1}\sum_{j=0}^{d_2} c_{i,j}(x,y) x^iy^j
\]
where the value of $c_{i,j}(x,y)$ depends only on the congruence classes 
of $x$ and $y$ modulo $m$.  Again, a quasi-polynomial in two variables, 
restricted to fixed congruence classes modulo $m$ for $x$ and $y$, 
is an ordinary polynomial.  
To compute the value of a quasi-polynomial $H(x,y)$ 
is as fast as computing the value of a polynomial:  
First determine the congruence classes of $x$ and $y$ modulo $m$, 
and then evaluate the appropriate polynomial.  

The following theorem  describes the construction of  quasi-polynomial packing functions on every integer sector.

\bt
Let $r$ and $s$ be relatively prime positive integers.
For ${\ell} \in \{0, 1, 2,\ldots, s-1 \}$, let 
\beq       \label{CantorSector:h1}
u_{\ell} = \left[ \frac{r{\ell}}{s} \right]
\eeq
and
\beq      \label{CantorSector:h2}
h^{({\ell})}_{r/s}(x,y) = \frac{r(x-{\ell})(x-{\ell}-s)}{2s^2} + \frac{(u_{\ell}+1)(x-{\ell})}{s} + y.
\eeq
The function $H_{r/s}(x,y)$ defined by 
\beq    \label{CantorSector:h3}
H_{r/s}(x,y) = sh^{({\ell})}_{r/s}(x,y) + {\ell}  \qquad \text{  if $x \equiv {\ell} \pmod{s}$ }
\eeq
is a quasi-polynomial packing function with period $s$ on the integral sector $I(r/s)$.  
\et

\begin{proof}
For $\ell \in \{ 0,1, 2,\ldots, s-1 \}$, let 
\[
\Omega_{\ell} = \{ (x,y) \in I(r/s) : x\equiv {\ell} \pmod{s} \}.
\] 
Then $\Omega_{\ell} \cap \Omega_{\ell'} = \emptyset$ if $\ell \neq \ell'$, 
and $\bigcup_{\ell=0}^{s-1} = I(r/s)$.  
We shall construct a packing polynomial $h^{({\ell})}_{r/s}(x,y)$ 
on each of the sets  $\Omega_{\ell} $, 
and then dilate and translate these polynomials to obtain 
a quasi-polynomial $H_{r/s}(x,y)$ with period $s$ and domain $I(r/s)$.  

Fix an integer $\ell \in \{ 0,1, 2,\ldots, s-1 \}$.  
For every nonnegative integer $a$, let 
\[
J_a = \{(sa+{\ell},y) : y = 0,1,\ldots, ra+u_{\ell}\}.
\]
Then $\bigcup_{a\in \N_0} J_a = \Omega_{\ell}$ and $J_a \cap J_{a'}$ if $a \neq a'$.  
We construct the polynomial $h^{({\ell})}_{r/s}(x,y)$ by enumerating 
the lattice points in $\Omega_{\ell}$ 
as follows:  
For $(\ell, y) \in J_0$, let $h^{({\ell})}_{r/s}(\ell,y) = y$ for $y = 0,1,\ldots, u_{\ell}$.
For $a \geq 1$, count the elements of $J_a$ after the elements of $J_{a-1}$, 
and count the elements of $J_a$ from the bottom lattice point $(sa+{\ell},0)$ 
to the top lattice point $(sa+{\ell},ra+u_{\ell})$.
If $x = sx'+\ell$, then 
\begin{align*}
h^{({\ell})}_{r/s}(x,y) & =h^{({\ell})}_{r/s}(sx'+\ell,y) \\
& =  \sum_{a=0}^{x' -1} |J_a| + y \\
& =  \sum_{a=0}^{x' -1} (ra+ u_{\ell} +1) + y  \\
& = \frac{rx'(x'-1)}{2} + (u_{\ell} +1) x' + y \\
& = \frac{r(x-{\ell})(x-{\ell}-s)}{2s^2} + \frac{(u_{\ell}+1)(x-{\ell})}{s} + y.
\end{align*}
The polynomial $h^{({\ell})}_{r/s}(x,y)$ is a packing polynomial on the set $\Omega_{\ell}$, 
and the quasi-polynomial $H_{r/s}(x,y) $ defined by~\eqref{CantorSector:h3} 
is a packing function on the integral sector $I(r/s)$.  
This completes the construction.
\end{proof}

For example, if $r=3$ and $s = 2$, then 
\[
h^{(0)}_{3/2}(x,y) = \frac{3x^2}{8} - \frac{x}{4} + y
\]
and
\[
h^{(1)}_{3/2}(x,y) = \frac{3x^2}{8} - \frac{x}{2} + y +  \frac{1}{8}.
\]
The quasi-polynomial 
\[
H_{3/2}(x,y) = 
\begin{cases}
 \frac{3x^2}{4} - \frac{x}{2} + 2y  & \text{if $x\equiv 0 \pmod{2}$ } \\
\\
 \frac{3x^2}{4} - x + 2y +  \frac{5}{4} & \text{if $x\equiv 1 \pmod{2}$. }
\end{cases}
\]
is a packing function on the integer sector $I(3/2)$.

\section{Open problems}
\benum
\item
Is there a quadratic packing polynomial for the sector $I(3/5)$?
Is there a quadratic packing polynomial for the sector $I(3/2)$?
For what rational numbers $\alpha$ do there exist quadratic packing polynomials?

\item
In Theorem~\ref{CantorSector:theorem:I(r)} we constructed two quadratic 
packing polynomials on the integer sector $I(r)$ for every integer $r \geq 2$?
Are these the only quadratic packing polynomials on $I(r)$

\item
Can any rational sector have more than two quadratic  packing polynomials?

\item
Can any rational sector have more than two packing  polynomials?

\item
Is there a rational sector with a packing polynomial  of degree greater than two?  

\item
Prove that  there is no  packing polynomial for any irrational sector. 

\eenum

\emph{Acknowledgements:}
I wish to thank for Dick Bumby and Tim Susse for helpful discussions.  

\def\cprime{$'$} \def\cprime{$'$} \def\cprime{$'$}
\providecommand{\bysame}{\leavevmode\hbox to3em{\hrulefill}\thinspace}
\providecommand{\MR}{\relax\ifhmode\unskip\space\fi MR }
\providecommand{\MRhref}[2]{%
  \href{http://www.ams.org/mathscinet-getitem?mr=#1}{#2}
}
\providecommand{\href}[2]{#2}

\end{document}